\documentclass[11pt]{article}
\usepackage{theorem,amssymb}
\advance \topmargin by -\headheight\advance \topmargin by -\headsep\textheight 8.9in\oddsidemargin 0pt\evensidemargin \oddsidemargin\marginparwidth 0.5in\textwidth 6.5in
\date{}
\newtheorem{theorem}{Theorem}
\newtheorem{lemma}{Lemma}
\newtheorem{corollary}{Corollary}

\newtheorem{definition}{Definition}
\theorembodyfont{\rm}
\newtheorem{remark}{Remark}
\begin{document} 
\begin{center}
\bf{\Large Integer-valued fixed point index for compositions of acyclic maps}
\vskip 15ptE. G. Sklyarenko,  G. S. Skordev 
\end{center}

\vskip 15pt
 
\begin{center}
\bf{\it Dedicated to the memory of Jean Leray}
\end{center}

\bigskip

\begin{abstract} 
An integer-valued fixed point index for compositions of  acyclic multivalued maps     is constructed. This integer-valued fixed point index has the properties: additivity, homotopy invariance, normalization, commutativity, multiplicativity. The  acyclicity is with respect to the \v{C}ech   cohomology   with integer coefficients. The technique of chain approximation is used.
\end{abstract} 

 \bigskip 
 
 {\it Keywords}: Fixed point index, multivalued maps, acyclic maps
 
\bigskip

The fixed point index (f. p. i.) for compositions of ${\mathbb F}$-acyclic multivalued maps of locally compact $ANR$ spaces was constructed in \cite{Sk1, Sk2, Sk3, Sk4, Dz, Lech}. This f. p. i. has values in the field ${\mathbb F}$ and has the properties: additivity, homotopy invariance, normalization, commutativity, multiplicativity and mod $p$. The ${\mathbb F}$-acyclicity of the multivalued maps is with respect to the \v{C}ech  homology, i.e., the reduced \v{C}ech  homology with coefficients in ${\mathbb F}$ of the images of all points is equal to zero.   The construction used the chain approximation technique. 

An integer-valued f. p. i. for ${\mathbb Z}$-acyclic multivalued maps of ENR spaces, i.e., finite dimensional ANR spaces, with all properties mentioned above, was constructed in \cite{Biel}, see also \cite{Granas}, p. 547-550. The ${\mathbb Z}$-acyclicity is with respect to the   \v{C}ech  cohomology. The construction of this integer-valued f. p. i. is based on some homotopic considerations. Both constructions are presented in \cite{Lech}, p. 163-173 and p. 251-276, where all necessary definitions and results used in this note are given, see also \cite{Granas}. 

The question was raised about the existence of integer valued f. p. i. with all properties for 
${\mathbb Z}$-acyclic maps, with respect to \v{C}ech  (co)homology, of compact ANR spaces, see \cite{Biel}, \cite{Lech}, ch. IV, Section 34, p. 173. The difficulty   is that the chain approximation technique, developed in \cite{Sk1, Sk2, Sk3, Sk4} and based on \cite{Begle1, Vietoris}, does not work in the case of ${\mathbb Z}$-acyclic maps with respect to the  \v{C}ech  homology, see the example in \cite{Begle2}. Moreover, the technique of \cite{Biel} is based on some results and arguments from   geometric topology for which   finite dimensionality of the spaces is essential.
A particular solution, based on the homotopy method,  was proposed in \cite{WK}. 

In this note we shall give a construction of integer-valued f. p. i. with all properties for   compositions of ${\mathbb Z}$-acyclic maps of compact ANR spaces. The acyclicity is with respect to the \v{C}ech  cohomology   with integer coefficients. The chain approximation technique is applied. For this a relation between \v{C}ech  cohomology theory with integer coefficients and the Steenrod-Sitnikov homology with integer coefficients and some results about the last homologies are used.

\section{Preliminaries}

Here we shall define some notions and fix some notations.
As usual, we denote by ${\mathbb N}$ the natural numbers and by ${\mathbb Z}$ the integers.

All topological spaces are assumed to be compact metric spaces. 

For a given space $X$ we denote by $Cov(X)$ the set of all its finite open coverings. For $\alpha , \beta \in Cov(X)$ we say that $\beta$ is a refinement of $\alpha$, denoted as $\beta > \alpha$, if for every element $U  \in \beta$ there is an element $V \in \alpha$ such that
$U \subset V$. 

Let $\alpha = \{U_0, \ldots , U_n\}\in Cov(X)$ and let $A$ be a closed subset of $X$. By 
$\alpha \mid A$ we denote the open covering of the set $A$ whose elements are all nonempty sets 
$U_0 \cap A , \ldots , U_n \cap A $.

If $A$ is a subset of $X$ and $\alpha \in Cov(X)$, we denote by $St(A, \alpha)$ the union of all elements of the covering $\alpha$ which meet $A$.   

We shall use the standard definitions for chain complexes, chain maps, chain homotopies, homologies, see e.g. \cite{Eilenberg, Spanier}. 
For compact  metric spaces we use \v{C}ech  cohomology and \v{C}ech   homology with coefficients in given Abelian groups, \cite{Eilenberg, Spanier, EG2}. We use also the Steenrod-Sitnikov homology theory with integer coefficients for compact metric spaces, \cite{EG1, EG2}.

For a given space $X$ and $\alpha \in Cov(X)$  we denote by $N(\alpha)$ the nerve of the covering $\alpha$. The vertices of this (abstract) simplicial complex are the elements of the covering $\alpha$. Furthermore, the $k+1$ vertices $U_0, \ldots , U_k$ are vertices of a $k$-simplex $\sigma ^k = [U_0, \ldots , U_k]$ in $N(\alpha)$
if and only if $U_0 \cap \ldots \cap U_k$ is not empty. The support $supp (\sigma ^k)$ of the simplex $\sigma^k$ is the set $U_0 \cup \ldots \cup U_k$. Let $c = \sum _{i=1}^l g_i \sigma _i^k \in C_k(N(\alpha); G), \, g_i \neq 0,$ be a chain of the simplicial complex $N(\alpha)$ with coefficients in an Abelian group $G$. Then the support $supp(c)$ of $c$ is the set
$ supp(\sigma_1^k) \cup \ldots \cup supp(\sigma_l^k)$. 

By $N(\alpha)^{(n)}$ we denote the $n$-th skeleton of the simplicial complex $N(\alpha)$.

Denote by $\pi (\beta , \alpha): N(\beta)\longrightarrow N(\alpha)$ an induced simplicial  map for   given coverings $\alpha ,\beta \in Cov(X)$ with $\beta > \alpha$. For a given closed subset $A$ in $X$ we denote also by 
$\pi (\beta , \alpha): N(\beta \mid A)\longrightarrow N(\alpha \mid A)$ the restriction of the map $\pi (\beta , \alpha)$ on the simplicial complex $N(\beta \mid A)$.

A multivalued map $F: X \longrightarrow Y$ is a map which assigns to every point $x \in X$ a nonempty compact  set $F(x)$ in $Y$. The graph $\Gamma (F)$ of the map $F$ is the set $\Gamma = \{(x,y) \in X \times Y : \, y \in F(x)\}$.
The map $F$ is called upper semi-continuous (u.s.c.) if the graph $\Gamma (F)$ is a closed subset in the space $X \times Y$. We have two projections $p: \Gamma (F)\longrightarrow X,$ defined by $p(x,y)=x$   and 
$q: \Gamma (F)\longrightarrow Y,$ defined by $q(x,y)=y$ for $(x,y) \in \Gamma(F)$. Then $F(x)= q (p^{-1}(x))$ for $x \in X$.
 For u.s.c. maps and their general properties see \cite{Berge, Lech}.

Let $\tilde{H}$ be a reduced homology or cohomology theory with coefficients in an Abelian group $G$. A compact space $Z$ is called $G$-acyclic with  respect to (w.r.t.) $\tilde{H}$ if $\tilde{H}(Z; G) =0$. An u.s.c. map $F: X \longrightarrow Y$ is called $G$-acyclic w.r.t. to $\tilde{H}$ if the spaces $F(x)$ are $G$-acyclic w.r.t. $\tilde{H}$ for all $x \in X$.

A single-valued continuous map $f: X \longrightarrow Y$ is called $G$-acyclic w.r.t. $\tilde{H}$ if the  inverse map
$f^{-1} : Y \longrightarrow X$ of $f$        is u.s.c. multivalued $G$-acyclic  w.r.t. $\tilde{H}$.  

Let $\Gamma (F)$ be the graph of the    $G$-acyclic w.r.t. $\tilde{H}$  multivalued  map $F$ and $p: \Gamma (F) \longrightarrow X,  \, q : \Gamma \longrightarrow Y$ are the projections defined above. Then the 
single-valued map $p: \Gamma (F) \longrightarrow X$ is $G$-acyclic w.r.t. $\tilde{H}$. The map $p: \Gamma (F) \longrightarrow X$ induces a homomorphism   in the (co)homologies $\tilde{H}$. This homomorphism is an isomorphism in the cases:

 (a) $\tilde{H}^{\ast}$ are the reduced \v{C}ech cohomologies with coefficients in an arbitrary Abelian group $G$, \cite{Spanier}, Theorem 15, p. 344.
  
(b) $\check{H}_{\ast}$ are the reduced \v{C}ech  homologies with coefficients in an arbitrary field $G$, \cite{Begle1}, Theorem 2, p. 538.  

(c) $\tilde{H}_{\ast}$ are the  reduced Steenrod-Sitnikov homology with integer  coefficients   \cite{Yanko}, Theorem   3.2., p. 57.

In the cases (b) and (c) we did not state the most general  results, for them  see \cite{EG2}, Ch. 8, 6.2.

In the case (a) we define the homomorphism 
$ F^{\ast} :  H^{\ast}(Y; G)\longrightarrow H^{\ast}(X; G) $ by $F^{\ast} = (p^{\ast})^{-1}q^{\ast}$, where $p^{\ast}: H^{\ast}(X; G)\longrightarrow H^{\ast}(\Gamma (F); G) $ and 
$q^{\ast}: H^{\ast}(Y; G)\longrightarrow H^{\ast}(\Gamma (F); G) $ are the homomorphisms induced by the projections $p$ and $q$, respectively.

Similarly, in the cases (b) and (c) we define  $ F_{\ast} : \check{H}_{\ast}(Y; G)\longrightarrow \check{H}_{\ast}(X; G) $  by $F_{\ast}=q_{\ast} (p_{\ast})^{-1}$, where 
$p_{\ast}: \check{H}_{\ast}(\Gamma (F); G)\longrightarrow \check{H}_{\ast}(X; G) $ and 
$q_{\ast}: \check{H}^{\ast}(\Gamma(F); G)\longrightarrow \check{H}^{\ast}(Y; G) $ are the homomorphisms induced by the maps $p$ and $q$, respectively.

\section{Maps of order $n$ with respect to $G$}

Here we shall give, in an appropriate form, some definitions and results of E. Begle, \cite{Begle1}.

\begin{definition}Let $X$ be a compact metric space and let $G$ be an Abelian group. Let $n$ be a natural number. The compact $X$ is called $(n, G)$-compact if for every covering $\alpha \in Cov(X)$ there is a covering $\mu \in Cov(X)$ such that $\mu > \alpha$ and the homomorphism
$$
\pi (\mu , \alpha)_{k}: \tilde{H}_k(N(\mu);G) \longrightarrow \tilde{H}_k(N(\alpha);G)
$$
is the zero homomorphism for $0 \leq k \leq n$.
\end{definition}

\begin{remark}
Every $(n, G)$-compact space is $(n, G)$-acyclic with respect to the \v{C}ech  homology, i.e., the reduced homology $ \check{H}_k(X; G)= 0$ for $k=0, \ldots , n$. 

The converse is not true, e.g., the $2$-adic solenoid is ($1$, ${\mathbb Z}$)-acyclic with respect to the \v{C}ech  homology with integer coefficients, but is not 
 $(1,  {\mathbb Z} )$-compact space, see \cite{Eilenberg}, Ch. X, Exercise F. Both properties are equivalent   when $G$ is a field.
\end{remark}

\begin{definition}(cf. \cite{Begle1}, Section 3)
Let $F: X \longrightarrow Y$ be an u.s.c. multivalued   map of the compact metric space $X$ in the compact metric space $Y$. Let $G$ be an Abelian group. The map  $F$ is called   map of order $n$ w.r.t. the group $G$, written (n, G)-map, if $F(x)$ is $(n, G)$-compact for   each point $x \in X$
\end{definition}

In \cite{Begle1} this definition is given for maps $F = f^{-1}$, where $f: Y \longrightarrow X$ is a single-valued continuous onto map. E. Begle called $f$ Vietoris map of order $n$. We say that $f$ is a single-valued $(n,  G)$-map.

Recall that E. Begle proved the Vietoris theorem for single-valued onto $(n,  G)$-maps   $f: Y \longrightarrow X$. For such maps  the induced homomorphism   $f_k : \tilde{H}_k (Y; G) \longrightarrow \tilde{H}_k (X; G)$ of the reduced \v{C}ech  homologies is an isomorphism for $0 \leq k \leq n$, \cite{Begle1}, Section 3, Theorem 1.
This theorem is not true in the case when the map $f$ is $(n, {\mathbb Z})$-acyclic with respect to the \v{C}ech  homology, \cite{Begle2}.  

  E. Begle derived the Vietoris theorem from   the two lemmas below, in the case of the map $F=f^{-1}$. 

\begin{lemma}(cf. \cite{Begle1}, Section 4, Lemma 2) If  $F: X \longrightarrow Y$ is an u.s.c. multivalued 
$(n,   G)$-map, then for each covering $\alpha \in Cov(X)$ and each covering $\beta \in Cov(Y)$ there is a covering 
$\nu = \nu (\alpha , \beta) \in Cov(X)$, with $\nu > \alpha$, and a chain map $T(\nu , \alpha)= \{T(\nu , \beta)_k\}, \, k= 0, \ldots , n+1$
$$
T(\nu, \alpha)_k :    C_k(N(\nu)^{(n+1)}; G) \longrightarrow C_k(N(\beta)^{(n+1)}; G)) 
$$
such that
\begin{enumerate}
\item
for any $k$-simplex $\sigma \in N(\nu)$ there is a point $x(\sigma) \in X$ with 
\begin{enumerate}
\item $supp (\sigma)\subset St(x(\sigma), \alpha)$,
\item $supp(T(\nu , \beta)_k(\sigma))\subset St(F(x(\sigma)), \beta)$,
\end{enumerate}
\item $f T(\nu , \alpha)_k(\sigma)$ is chain homotopic to $\sigma$ on $N(\alpha)$.
\end{enumerate}
\end{lemma}

For the next Lemma we use the notations of Lemma 1.

\begin{lemma}(cf. \cite{Begle1}, Section 4, Lemma 3) Let  $F: X \longrightarrow Y$ be an u.s.c. multivalued  
$(n,   G)$-map. Let 
$\alpha$ and $ \alpha _1$ be coverings of $X$ with $\alpha _1 > \alpha$. Let $\beta$ and $\beta _1$ be coverings of $Y$ with $\beta _1 > \beta$. Let $\nu = \nu (\alpha , \beta)$ and $\nu _1= \nu (\alpha _1, \beta _1)$. Let $T(\nu , \beta)$ and 
$T_1(\nu _1, \beta _1)$ be the chain maps from Lemma 1. Then there is a common refinement $\gamma $ of the coverings $\nu$ and $\nu _1$ such that the chain map  $T(\nu , \beta ) \pi(\gamma , \nu)$ is chain homotopic with
$\pi (\beta _1 , \beta)T(\nu _1 , \beta _1 ) \pi(\gamma , \nu _1)$ to a chain homotopy $D$ with the property:
for every simplex $\sigma \in N(\gamma)$ there is a point $c(\sigma) \in X$ with:
\begin{enumerate}
\item $supp (\sigma) \subset St(c(\sigma), \alpha)$,
\item $supp (D(\sigma))\subset St(F(c(\sigma)), \beta)$.
\end{enumerate}
We say that the homotopy $D$ is $(F, \alpha , \beta )$-small.
\end{lemma}

\begin{remark}
1. Lemma 2 and Lemma 3 in \cite{Begle1} are stated for the multivalued map $F=f^{-1}  $, i.e., for the inverse of the   single-valued map $f: Y \longrightarrow X$ in an equivalent form for Vietoris chains.
The more general form given above  follows easily from   Lemma 2 and Lemma 3 in \cite{Begle1}.

2. The chain maps $\{ T(\nu, \beta), \alpha   \in Cov(X), \beta \in Cov(Y), \nu =  \nu ( \alpha , \beta) \in Cov (X) \}$ induce  the homomorphisms 
$F_k: \check{H}_k(X; G) \longrightarrow \check{H}_k(Y; G)$
for $0 \leq k \leq n$.

3. The properties 1a),b) of Lemma 1 and 1, 2 of Lemma 2 are not stated in Lemma 2, Lemma 3 in \cite{Begle1}, in the  case $F=f^{-1}$, but are explicit in  the proofs of these lemmas given there. 
\end{remark}

\section{Chain approximations and approximation systems for   u.s.c. maps}

Chain approximations for u.s.c. multivalued maps were used by L. Vietoris (see \cite{Vietoris}), S. Eilenberg and D. Montgomery (\cite{Eilenberg1}), E. Begle (\cite{Begle1, Begle3}) and B. O'Neil (\cite{Neil}). The explicit definition is given in \cite{Sk1}.  They are developed further in \cite{Sk2, Sk3, Sk4, Dz}, \cite{Lech}, Ch. 4, p. 251-276.
In all these papers the authors, except E. Begle in \cite{Begle1}, work  with $G$-acyclic maps w.r.t.  the \v{C}ech  homology with coefficients in a field $G$. The most general case is considered by E. Begle in \cite{Begle1} for \v{C}ech  homology with coefficients in an Abelian group $G$ for   single-valued
$(n, G)$-maps. As mentioned before every $(n,  G)$-map is $(n, G)$-acyclic but not vice versa.

 Here we shall give the definitions of chain approximations and approximation systems for $(n,  G)$-maps. They are the same as the definitions given in \cite{Dz}, \cite{Lech},   Ch. 4, p. 251-276, in the case where $G$ is a field.
 
The chain approximations and the approximation systems for $(n,   {\mathbb Z})$-maps are defined in the same way and have the same properties as the chain approximations and the approximation systems for 
${\mathbb F}$-acyclic maps
for ${\mathbb F}$ a field. Moreover, the proofs of the corresponding properties for $(n, G)$-maps are the same as for ${\mathbb F}$-acyclic maps. For this reason we skip these proofs, but give exact references for the corresponding proofs in \cite{Lech}.

\begin{definition}(cf.   \cite{Lech}, Definition 50.29, p. 255)
Let $F: X \longrightarrow Y$  be an u.s.c. multivalued map and let $G$ be an  Abelian group. Let 
$\alpha , \overline{\alpha} \in Cov(X)$, $\overline{\alpha} > \alpha$, and let $\beta \in Cov(Y)$. An augmentation preserving  chain map 
$\varphi : C_{\ast}(N(\overline{\alpha}) ^{(n+1)}; G) \longrightarrow C_{\ast}(N(\beta) ^{(n+1)}; G)$ is called 
$(n, \alpha , \beta)$-approximation of the map $F$ provided for each simplex $\sigma \in N(\overline{\alpha}) ^{(n+1)}$
there is a point $x(\sigma) \in X$ such that 
\begin{enumerate}
\item  
$supp (\sigma)\subset St(x(\sigma), \alpha)$,
\item 
$supp (\varphi (\sigma)) \subset St(F(x(\sigma)), \beta)$.
\end{enumerate}
\end{definition}

\begin{remark}
Let $F: X \longrightarrow Y$ be an u.s.c. multivalued $(n,   G)$-map. 
The chain map $T(\nu , \beta)$ from Lemma 1 is an $(n, \alpha , \beta )$-approximation of the map 
$F$.  
\end{remark}

\begin{definition}
Let $F: X \longrightarrow Y$ be an u.s.c. multivalued map. An $(n,G)$-approximation system, written $(n, G, A)$-system, of the map $F$ is a collection of chain maps 
$\{ \varphi (\nu , \beta) : \, \,  \alpha , \nu=\nu(\alpha , \beta) \in Cov (X), \beta \in Cov(Y)\}$, where
\begin{itemize}
\item the chain map $\varphi (\nu , \beta) : C_{\ast}(N(\nu)^{(n+1)}; G) \longrightarrow C_{\ast}(N(\beta)^{(n+1)}; G)$ is
an    $(n, \alpha , \beta)$-approximation of $F$. 
\item  Let  $\nu = \nu(\alpha , \beta) \in Cov(X)$  correspond to a given  $\alpha \in Cov(X)$ and $\beta \in Cov(Y)$. Let $\nu _1= \nu(\alpha _1, \beta _1) \in Cov(X)$   correspond to a given $\alpha _1 > \alpha$ and  
$  \beta _1 > \beta  $. Then it follows: 
there exists $\gamma \in Cov(X)$ wit $\gamma > \nu $ and $\gamma > \nu _1$ such that the chain maps
$ \varphi (\nu , \beta) \pi (\gamma , \nu)$ and 
$\pi (\beta _1 , \beta)\varphi (\nu _1, \beta _1)\pi (\gamma , \beta _1)$ are homotopic with a chain homotopy which is $(F, \alpha, \beta)$-small.
\end{itemize}
\end{definition}

\begin{remark}

1. Let $F: X \longrightarrow Y$ be an u.s.c. multivalued $(n,   G)$-map. 
Lemma 1 and Lemma 2 imply that the collection of chain maps 
${\cal A}(F)= \{T(\nu , \beta) : \, \, \alpha \in Cov(X), \, \nu = \nu (\alpha  , \,  \beta) \in Cov(X), \beta \in Cov(Y)\}$ is an 
$(n, G, A)$-system for the map F. We call ${\cal A}(F)$ the induced $(n, G, A)$-system of the map $F$.

2. For the definition of an $(n,  G, A)$-system it is not necessary to consider all coverings of the spaces $X$ and $Y$. It is enough to consider only a fundamental sequence of coverings $\{ \alpha _k \}  \in Cov(X)$ and 
$\{ \beta _k \}     \in Cov(Y)$. Then an $(n, G, A)$-system is a collection of augmentation preserving chain maps
$\{ \varphi (l, k)= \varphi(\nu _l, \beta _k): k,  l \in {\mathbb N} \}$, which satisfy the conditions of the previous definition with $\nu _l = \nu (\alpha _k , \beta _k)$. In the case where $X$ and $Y$ are finite polyhedra with  given triangulations $\tau$, $\mu$, respectively, we take the covering  $\alpha _k$ to consist of all open stars of the vertices of the $k$-th barycentric subdivision $\tau ^k$ of the triangulation $\tau$ w.r.t. $\tau ^k$. Similarly for $\beta _k$. Compare with     \cite{Lech}, Definitions 50.17, 50.18, p. 152.
\end{remark}

Now we shall consider a composition of u.s.c. multivalued maps and shall define a composition of $(n, G, A)$-systems. 
From   \cite{Lech}, p. 259-262 follows:

\begin{lemma}(cf. \cite{Lech}, Proposition (50.37), p. 260) Let $F_i : X_i \longrightarrow X_{i+1}, \, i = 1, 2,$ be u.s.c. multivalued maps. Let  
$\alpha \in Cov(X_1), \,  \gamma \in Cov(X_3)$. There exists $\beta \in Cov(X_2)$ such that if 
$\varphi (\nu _2 , \gamma) : C_{\ast}(N(\nu _2)^{(n+1)}; G) \longrightarrow C_{\ast}(N(\gamma)^{(n+1)}; G)$ is an 
$(n, \beta , \gamma)$-approximation of $F_2$ for $\nu _2 = \nu (\beta ,  \gamma) > \beta$, and 
$\varphi (\nu _1 , \nu _2) : C_{\ast}(N(\nu _1)^{(n+1)}; G) \longrightarrow C_{\ast}(N(\nu _2)^{(n+1)}; G)$ is an 
$(n, \alpha , \nu _2)$-approximation of $F_2$ for $\nu _1 = \nu (\alpha ,  \nu _2) > \alpha$, then
$ \varphi(\nu _2 , \gamma) \varphi (\nu _1, \nu _2) : C_{\ast}(N(\nu _1)^{(n+1)}; G)\longrightarrow 
C_{\ast}(N( \gamma  )^{(n+1)}; G)$ is an $(n, \alpha , \gamma )$-approximation of the map $F_2F_1$.
\end{lemma}

\begin{lemma}(cf. \cite{Lech}, Proposition (50. 39), p. 261) Let $F_i : X_i \longrightarrow X_{i+1}, i = 1, 2,$ be u.s.c. multivalued maps.  Assume that
${\cal A}_1 = \{ \varphi (\nu _1 , \beta) : \, \, \alpha \in Cov(X_1), \nu _1= \nu(\alpha , \beta) \in Cov(X_1), \beta \in Cov(X_2)  \}$ is an $(n, G, A)$-system of the map $F_1$
and 
${\cal A}_2 = \{ \varphi (\nu _2 , \gamma) : \, \, \beta \in Cov(X_2), \nu _2= \nu(\beta , \gamma) \in Cov(X_2), \gamma \in Cov(X_3)  \}$   is an $(n, G, A)$-system of the map $F_2$. Then 
$\{ \varphi(\nu _2 , \gamma) \varphi (\nu _1, \nu _2): \, \alpha \in Cov(X_1), \gamma \in Cov(X_3)\}$
is an $(n, G, A)$-system of the map $F_2F_1$.
\end{lemma}

We call this $(n, G, A)$-system for the map $F_2F_1$   the composition of the 
${\cal A}$-systems ${\cal A}_1$ and ${\cal A}_2$ and denote it by ${\cal A}_2 \circ {\cal A}_1$
 
From this lemmas and Remark 3, 4 follows

\begin{corollary} Let $F_i: X_i \longrightarrow X_{i+1}, \, i= 1, \ldots , k-1$ be $(n,   G)$-maps. Let 
${\cal A}(F_i)$ be the induced $(n, G, A)$-system of the map $F_i$, defined in Remark 4.1. Then the composition
${\cal A}(F_{k-1}) \circ \ldots \circ {\cal A}(F_1)$ is an $(n, G, A)$-system of the map 
$F=F_{k-1}\ldots F_1 : X_1 \longrightarrow X_k$.
\end{corollary}
We denote this $(n, G,   A)$-system of $F= F_{k-1}\ldots F_1$     by ${\cal A}(F)$.

\begin{definition}(cf. \cite{Lech}, Definition 50.19, p. 252) Let $F_1, F_2 : X \longrightarrow Y$ be u.s.c. multivalued maps. Let $H: X \times I \longrightarrow Y$ be an u.s.c. multivalued homotopy joining $F_1$ and $F_2$. Let ${\cal A}(F_i) $ be an $(n, G, A)$-system   of $F_i, \, i=1,2$. The approximation  systems ${\cal A}(F_1)$ and ${\cal A}(F_2)$ are called $H$-homotopic if for all sufficiently fine coverings $\alpha \in Cov(X), \, \beta \in Cov(Y)$ and
$\varphi_1(\nu _1 , \beta) \in {\cal A}(F_1), \, \varphi_2(\nu _2 , \beta) \in {\cal A}(F_2)$ the chain maps 
$\varphi_1(\nu _1 , \beta) \pi(\gamma , \nu _1)  $ and $\varphi_2(\nu _2 , \beta) \pi(\gamma , \nu _2)$ are chain homotopic with a chain homotopy $D$ for all $\gamma > \nu _1 , \nu _2$. Moreover, we assume that the chain homotopy $D$ satisfies the following condition: 
for every simplex $\sigma \in N(\gamma)$ there is a point $d(\sigma) \in X$ such that
\begin{itemize}
\item $supp (\sigma) \subset  St(d(\sigma), \alpha)$, 
\item $supp (D \sigma) \subset St(H(d(\sigma) \times I), \beta)$.
\end{itemize}
\end{definition}

\begin{lemma}(cf. \cite{Lech}, Lemma 51.8, p. 265) Let $F_1, F_2 : X \longrightarrow Y$ be $(n,   G)$-maps. Let $H: X \times I \longrightarrow Y$ be an $(n,   G)$-map, which is a  homotopy joining $F_1$ and $F_2$. Let ${\cal A}(F_i)$ be the induced $(n, G, A)$- system   of $F_i, \, i=1,2$. Then the  approximation  systems ${\cal A}(F_1)$ and ${\cal A}(F_2)$ are  $H$-homotopic.
\end{lemma}

\begin{definition}Assume that $F = F_{k-1} \ldots F_1$, and $\Phi = \Phi _{k-1}\ldots \Phi _1 : X_1 \longrightarrow X_{k}$ are compositions of the maps $F_i, \Phi _i : X_i \longrightarrow  X_{i+1}, \, i=1, \ldots , k-1$. Assume that 
$F_i$ and $\Phi _i$ are homotopic with a homotopy $H_i : X_i \times I \longrightarrow X_{i+1} , \, i= 1, \ldots , k-1$. Then we say that the maps $F$ and $\Phi$ are composition-homotopic  with a homotopy     
$H = H_{k-1}(H_{k-2} \times id)   \ldots (H_1 \times id)$, where $id: I \longrightarrow I$ is the identity.
\end{definition}

\begin{corollary}(cf. \cite{Lech}, Proposition 51.9, p. 265)
Assume that $F = F_{k-1} \ldots F_1, \Phi = \Phi _{k-1}\ldots \Phi _1 : X_1 \longrightarrow X_{k}$ are compositions of the $(n,   G)$-maps $F_i, \Phi _i : X_i \longrightarrow  X_{i+1}, \, i=1, \ldots , k-1$. Assume that 
$F$ and $\Phi $ are composition-homotopic  with a homotopy $H$ Then the  induced $(n, G, A)$-systems ${\cal A}(F)= {\cal A}(F)_{k-1}\circ \ldots \circ {\cal A}(F)_1$ and ${\cal A}(\Phi)= {\cal A}(\Phi)_{k-1}\circ \ldots \circ {\cal A}(\Phi)_1$ are 
$H$-homotopic.
\end{corollary}

\section{Fixed point index of an approximation system}

Let $(K, \tau )$ be a finite simplicial complex with triangulation $\tau$. By $\tau ^k$ we denote the $k$-th barycentric subdivision of the triangulation $\tau$.   Let 
$C_{\ast}(\tau ^k)= C_{\ast}(K, \tau ^k ; {\mathbb Z})$ be the chain complex with integer coefficients of the triangulation $\tau ^k$. By 
$b(k, l): C_{\ast}(\tau ^k)\longrightarrow C_{\ast}(\tau ^l)$ we denote the barycentric subdivision operation, $l >k$. We consider the fundamental sequence of open coverings $\{ \alpha _k \} $ of the space $K$, where $\alpha _k$ is the covering of $K$ by the open stars of the vertices of $\tau ^k$, w.r.t. the triangulation $\tau ^k$.

Let $U$ be an open set in $K$ such that its closure $\overline{U}$ is a subcomplex in $(K, \tau)$. By
$p(U, k) = \{p(U, k)_i\} : C_{\ast}(\tau ^k)\longrightarrow C_{\ast}(\tau ^k \mid \overline{U})$ we denote the projection homomorphism. Here $\tau ^k \mid \overline{U}$ is the restriction of the triangulation $\tau ^k$ on $\overline{U}$.

Let $F: \overline{U}\longrightarrow K$ be an u.s.c. multivalued map without fixed points on the boundary $\partial U$ of the set $U$, i.e., $x \notin F(x)$ for $x \in \partial U$, or the same, $Fix(F)= \{x : \, x \in F(x)\} \subset U$. In this case we call the triple $(K, F, U)$ admissible.

\begin{definition}We call a quadruple $(K, F, U; {\cal A})$  admissible if  $(K, F, U)$ is an admissible triple and ${\cal A}$ is an $(n, {\mathbb Z}, A)$-system of the map $F$, with $n =dim \, K$. Denote by ${\cal K}_{\cal A}$ the set of all  admissible quadruples.
\end{definition}

\begin{definition}(Fixed  point index on ${\cal K}_{\cal A}$, cf. \cite{Lech}, Ch. 4, Definition 50.21) The fixed point index
$I_{\cal A}: {\cal K}_{\cal A} \longrightarrow {\mathbb Z}$ is defined as follows. Let 
$(K, F, U; {\cal A}) \in {\cal K}_{\cal A}$ and let 
${\cal A}= \{ \varphi (l, k), \nu _l = \nu(\alpha _k , \alpha _k)\}$ be an $(n, {\mathbb Z}, A)$-system of the map 
$F: \overline{U} \longrightarrow K$, see Remark 4.2. Let  $\psi$ be the graded homomorphism
$\psi  = \{\psi _i\} = p(U, k) \varphi (l, k) b(k,l): C_{\ast}(\tau ^k \mid \overline{U}) \longrightarrow C_{\ast}(\tau ^k \mid \overline{U})$ and let $k$ be a sufficiently large natural number. 
  Then   $I_{{\cal A}}(K, F, U)$ is defined   by the Lefschetz number of the graded homomorphism $\psi$, i.e., 
$I_{{\cal A}}(K, F, U) = \lambda (\psi )= \sum _{i} (-1)^i tr(\psi _i)$,   where $tr(\psi  _i)$ is the trace of the homomorphism $\psi _i$.
\end{definition}
 
 \begin{remark}
 The definition of the fixed point index $I_{{\cal A}}(K, F, U)$ is correct, i.e., it does not depend on the number $k$. This follows as in Lemma and Definition 1.2, \cite{Sk1}.
\end{remark}

\begin{lemma} The fixed point index $I_{{\cal A}}(K, F, U)$ has the following properties:
\begin{enumerate}
\item {\bf Additivity}

Let $U_1, U_2$ be open, disjoint and polyhedral subsets of $U$ and   $Fix(F)= \{x : \, x \in F(x)\}\subset U_1 \cup U_2$. Then
$$
 I_{{\cal A}}(K, F, U) = I_{{\cal A}}(K, F, U_1) + I_{{\cal A}}(K, F, U_2).
$$

\item {\bf Homotopy invariance }

Let $H = H_t: \overline{U} \times I \longrightarrow K$ be an u.s.c homotopy  such that $(K, H_t, U)$ is an admissible triple for all $t \in I$. Let ${\cal A}_0, {\cal A}_1$ be $H$-homotopic 
$(n, {\mathbb Z}, A)$-systems for the maps $H_0, H_1$, respectively. Then 
$$
I_{{\cal A}_{0}}(K, F, U) = I_{{\cal A}_{1}}(K, F, U).
$$

\item {\bf Commutativity} 

Let $K, L$ be finite simplicial complexes. Let $W \subset K$ be an open subset  and let $F_1: K \longrightarrow L$,
$F_2: L \longrightarrow K$ be u.s.c. multivalued maps. Assume that $x \notin F_2F_1(x)$ for 
$x \in \partial{W}$ and 
$y \notin F_1F_2(y)$ for $y \in \partial{F_2^{-1}(W)}$. Assume further   that  
$$
y \in Fix(F_1F_2) \setminus \overline{F_2^{-1}(W)} \, \,  \mbox{implies} \, \,   F_2(y) \cap Fix(F_2 F_1 \mid \overline{W}) = \emptyset.
$$
Then for all $(n, {\mathbb Z}, A)$-systems ${\cal A}_1$,  ${\cal A}_2$ of $F_1, F_2$, respectively, follows
$$
I_{{\cal A}_{1} \circ {\cal A}_{2}}(L, F_1F_2, F_2^{-1}(W)) = I_{{\cal A}_{2} \circ {\cal A}_{1}}(K, F_2F_1, W).
$$
Here $n \geq \max \{dim \, K , dim \, L\}$.

\item {\bf Normalization }

$$
I_{\cal A}(K, F, K) = \lambda (\psi _{\ast}),
$$
where $\psi _{\ast} : H_{\ast}(K; {\mathbb Z}) \longrightarrow H_{\ast}(K; {\mathbb Z})$ is the homomorphism induced by the chain map $\varphi (l, k) b(k, l)$, where $\varphi (k, l) \in {\cal A}$ and $k$
is sufficiently large.
\end{enumerate}
\end{lemma}

Now we define the fixed point index $I_{\cal A}(K, F, V)$ for all open sets $V$ and maps 
$F: \overline{V} \longrightarrow K$ such that
  $Fix(F) \subset V$. We call such a triple admissible. 

Take an open polyhedral set $U$, such that $\overline{U} \subset V $ and $Fix(F) \cap (V \setminus U) = \emptyset$. Then define $I_{\cal A}(K, F, V) = I_{\cal A}(K, F, U)$. 

The proofs of Propositions  1.5, 1.6, \cite{Sk1} imply the Homotopy invariance and the Commutativity property.

The Normalization follows from \cite{Hopf}. 

\begin{definition}
Let ${\cal K}(n) = \{(K, F, U)\}$ where 
\begin{itemize}
\item $K$ is a finite simplicial complex, and the triple $(K, F, U)$ is admissible;
\item  the map $F= F_{k-1}\ldots F_1$ with $F_i :X_i \longrightarrow X_{i+1}$, $X_1 = \overline{U}$, 
$X_{k}=K$;
\item the maps $F_i$ are $(n,   {\mathbb Z})$-maps, $n = dim \,  K$, see Definition 2. 
\end{itemize} 
\end{definition}

\begin{definition}Let $(K, F, U) \in {\cal K}(n)$. Let 
${\cal A}(F) = {\cal A}(F_{k-1})\circ \ldots \circ {\cal A}(F_1)$ be the induced $(n, G, A)$-system of the map $F$. Then $(K, F, U; {\cal A}(F)) \in {\cal K}_{\cal A}$, see Definition 7. The fixed point index $I : {\cal K}(n) \longrightarrow {\mathbb Z}$ is defined as
$I(K, F, U) = I_{{\cal A}(F)}(K, F, U)$.
\end{definition}

From Lemma  6   follows 

\begin{corollary}
The fixed point index $I : {\cal K}(n) \longrightarrow {\mathbb Z}$ has the properties Additivity, Homotopy invariance, Commutativity and Normalization.
\end{corollary}

\begin{remark}
For the definition of the fixed point index $I: {\cal K}(n) \longrightarrow {\mathbb Z}$   a block structure of the simplicial complex could be used. Then one   obtains also the Multiplicativity of $I(K, F, U)$, as in  \cite{Sk4}.
\end{remark}

\section{${\mathbb Z}$-acyclic and $(n, G)$-maps}

Here we shall describe   relations between the ${\mathbb Z}$-acyclicity w.r.t. the  \v{C}ech  cohomology, w.r.t.  the Steenrod-Sitnikov homology   of compact metric spaces and $(n, {\mathbb Z})$-compact spaces, see Definition 1.

\begin{lemma} A compact metric space is ${\mathbb Z}$-acyclic  w.r.t. the \v{C}ech  cohomology      if and only if it is ${\mathbb Z}$-acyclic w.r.t. the Steenrod-Sitnikov homology.
\end{lemma}
Proof. Assume that the compact metric space is ${\mathbb Z}$-acyclic w.r.t. the \v{C}ech  cohomology. Then the exact sequence 
$$
0 \longrightarrow Ext(H^{n+1}(X, {\mathbb Z}); {\mathbb Z}) \longrightarrow H_n(X; {\mathbb Z}) \longrightarrow Hom(H^n(X; {\mathbb Z}), {\mathbb Z})\longrightarrow 0, 
$$ 
see \cite{EG2}, (14), p. 221, implies that $X$ is ${\mathbb Z}$-acyclic w.r.t. the   Steenrod-Sitnikov homology.  

For the inverse assertion assume that the compact $X$ is ${\mathbb Z}$-acyclic w.r.t. the   Steenrod-Sitnikov homology. The above exact sequence gives that $Ext(H^{n+1}(X, {\mathbb Z}); {\mathbb Z})=0$ and
$Hom(H^n(X; {\mathbb Z}), {\mathbb Z})=0$. Since $X$ is a compact metric space the groups 
$H^n(X; {\mathbb Z})$ are countable. Then the theorem of Stein-Serre implies that the groups 
$H^n(X; {\mathbb Z})$ are free, see \cite{EG3}, Proposition 1.2, Proposition 1.3 and Remark on p. 374. Then $Hom(H^n(X; {\mathbb Z}), {\mathbb Z})=0$ gives that the space $X$ is ${\mathbb Z}$-acyclic w.r.t. the \v{C}ech  cohomology.

\begin{lemma} Let $F: X \longrightarrow Y$ be an u.s.c. multivalued map. Assume that $F$ is $G$-acyclic w.r.t. the Steenrod-Sitnikov homology   with   coefficients in a countable Abelian group $G$.  Then $F$ is an
$(n,  G)$-map for every natural number $n$.
\end{lemma}
Proof.
Consider the representation  $F= q p^{-1}$ with the projections $p: \Gamma (F) \longrightarrow X$ and 
$q: \Gamma (F) \longrightarrow Y$. The Proposition 1.5.  from \cite{Yanko} gives that $p^{-1}(x)$ are $(n,   G)$-compact spaces  for $x \in X$ and every $n$. Then by   Definition  2 the map $F$ is   an $(n,     G)$-map for every $n$.

\begin{remark}
In the proof of Proposition 1.5.  from \cite{Yanko}  is used that the Mittag-Leffler property of projective systems $\{ A_i\}$ of Abelian groups is equivalent to the vanishing of the first derived functor $lim^1\{ A_i\} $ of the projective  system $\{ A_i\}$. This is true for   countable projective systems of countable Abelian groups, see \cite{EG3}, Proposition 1.2, p. 371. For this reason we assume that the spaces are compact and metric.
\end{remark}

 \section{Main Theorem}

\begin{definition}
Let $X$ be a compact ANR and let $U $ be an open set in $X$.  The triple $(X, F, U)$ is called acyclic admissible if 
\begin{itemize}
\item $F: \overline{U} \longrightarrow X$ has no fixed points on the boundary $\partial U$ of the set $U$, i.e., $Fix(F) \subset  U  $;
\item There is a natural number $k$ such that $F=F_{k-1}\ldots F_1 : \overline{U}\longrightarrow X$ and 
$F_i : X_i \longrightarrow X_{i+1} , \, \, i= 1, \ldots , k-1$, $X_1 = \overline{U}, X_{k}=X$;
\item the maps $F_i, \, i= 1, \ldots , k-1$ are ${\mathbb Z}$-acyclic w.r.t. \v{C}ech  cohomology with integer coefficients.
\end{itemize}
Denote by ${\cal K}$ the set of all acyclic admissible triples.
\end{definition}
For $(X, F, U)\in {\cal K}$ we have a homomorphism $F^{\ast}= F^{\ast}_1 \ldots F^{\ast}_{k-1} : H^{\ast}(X; {\mathbb Z}) \longrightarrow H^{\ast}(\overline{U}; {\mathbb Z})$, where $ F^{\ast}_i : H^{\ast}(X_{i+1}; {\mathbb Z}) \longrightarrow H^{\ast}(X_i; {\mathbb Z})$  is defined as follows. Consider the representation $F_i = q_i (p_i)^{-1}$ of the map with the projections $p_i : \Gamma (F_i) \longrightarrow X_i$  and $q_i : \Gamma (F_i) \longrightarrow X_{i+1}$. The Vietoris theorem implies that the homomorphism $p^{\ast}_i : H^{\ast}(X_i; {\mathbb Z}) \longrightarrow H^{\ast}(\Gamma (F_i); {\mathbb Z})$
is an isomorphism, see \cite{Spanier}, Theorem 15, p. 344. Then 
$F_i^{\ast} = (p_i^{\ast})^{-1}q_i^{\ast}: H^{\ast}(X_{i+1}) \longrightarrow H^{\ast}(X_{i})$, and 
$F^{\ast}= F_1^{\ast} \ldots F_{k-1}^{\ast}: H^{\ast}(X; {\mathbb Z}) \longrightarrow H^{\ast}(\overline{U}; {\mathbb Z})$.

Consider the case $X=U$.
Since $X$ is compact ANR then $H^{\ast}(X; {\mathbb Z})$ is a finetely generated Abelian group. Then the Lefschetz number $\lambda (F^{\ast})= \sum  _k (-1)^k tr(F^{k})$ is defined.  Here $F^{\ast}=\{ F^{k} \}$ and 
$F^{k} : H^k(X; {\mathbb Z})\longrightarrow H^k(X; {\mathbb Z})$ is the   homomorphism induced by the  map $F$.

Let $(X, F, U), (X, G, U) \in {\cal K}$. We call the triples  $(X, F, U$) and $(X, G, U)$ admissible homotopic if the maps $F$ and $G$ are composition-homotopic  with a homotopy $H$ such that $x \notin H(x,  t)$ for all $x \in \partial U$ and $t \in I= [0, 1]$, see Definition 6.

\begin{theorem}(integer fixed point index on {\cal K}) There is a function   $i: {\cal K}\longrightarrow {\mathbb Z}$ defined on the set of all acyclic admissible triples ${\cal K}$ with integer values with the following properties:
\begin{enumerate} 
\item {\bf Additivity}

Let $(X, F, U) \in {\cal K}$ and let $U_1, U_2 $ be open subsets of $U$ with $Fix(F) \subset U_1 \cup U_2$, then
$$
i(X, F, U) = i(X, F, U_1)+ i(X, F, U_2);
$$
\item {\bf Homotopy invariance}

Let $(X, F, U), (X, G, U) \in {\cal K}$ be admissible homotopic triples, then
$$
i(X, F, U) = i(X, G, U);
$$
\item {\bf Commutativity}

Let $F: X \longrightarrow Y$, $G: Y \longrightarrow X$ be compositions of acyclic maps. Let $U$ be an open set in $X$. Assume that $(X, GF,  U), (Y, FG, G^{-1}(U)) \in {\cal K}$ and $G(Fix(FG) \setminus G^{-1}(U)) \cap Fix(GF \mid U) = \emptyset$, then
$$
i(X, GF, U) = i(Y, FG, G^{-1}(U));
$$

\item {\bf Normalization}

$$
i(X, F, X) = \lambda (F^{\ast}).
$$
\end{enumerate}
\end{theorem}
Proof.

Let ${\cal K}'$ be the subset of ${\cal K}$ consisting of all acyclic admissible triples $(K, F, U) \in {\cal K}$ with $K$ a  finite simplicial complex.

Lemma 8 implies that if $(K, F, U) \in {\cal K}'$ then $(K, F, U) \in {\cal K}(n)$ with $n = dim \, K$, see Definition 9.

Consider $I: {\cal K}(n) \longrightarrow {\mathbb Z}$, see  Definition 10. Define $i(K, F, U)= I(K, F, U)$. 
Lemma 6 implies that the fixed point index $i: {\cal K}' \longrightarrow {\mathbb Z}$ has Additivity, Homotopy invariance and Commutativity properties. 

The  Normalization property in Lemma 6 is stated as follows:
$i(K, F, K) = \lambda (F_{\ast})$, where $F_{\ast} = (F_{k-1})_{\ast} \ldots (F_1)_{\ast}$ is a  homomorphism in the \v{C}ech  homology with integer coefficients.
It follows that $\lambda (F_{\ast}) = \lambda (F^{\ast})$, i.e., the Normalization property from   Theorem 1 also follows for the index $i: {\cal K}' \longrightarrow {\mathbb Z}$.

Since the  fixed point index $i: {\cal K}' \longrightarrow {\mathbb Z}$ has the Commutativity property we  can apply the extension   procedure from \cite{Lech}, Ch. 4, Section 53. As a result we obtain a fixed point index $i: {\cal K} \longrightarrow {\mathbb Z}$ with Additivity, Homotopy invariance, Commutativity and Normalization properties.

\begin{remark}

1. The fixed point index $i(X, F, U)$ coincides with the fixed point index defined in \cite{Biel} for $dim \, X < \infty$ and $F$   a  ${\mathbb Z}$-acyclic map.

2. Using the technique from \cite{Sk4}  one can prove that the fixed point index
$i: {\cal K} \longrightarrow {\mathbb Z}$ has also the Multiplicativity property, i.e., for $(X_1, F_1, U_1), (X_2, F_2, U_2) \in {\cal K}$ follows $(X_1 \times X_2, F_1 \times F_2, U_1 \times U_2) \in {\cal K}$ and 
$$
i(X_1 \times X_2, F_1 \times F_2,  U_1 \times U_2) = i(X_1, F_1, U_1)i(X_2, F_2, U_2)
$$

3. With an adaptation of the technique of \cite{Fritz1, Fritz2} one can obtain an integer fixed point index for compositions of multivalued, ${\mathbb Z}$-acyclic, weighted maps with an integer multiplicity function. For these maps see \cite{Pej}.
\end{remark}

{\bf Acknowledgments:} The paper was written during the visit of the first author in Bremen. The authors thank Heinz-Otto Peitgen for his  support and Jean-Paul Allouche for many valuable remarks. Many thanks also to Lech Gorniewicz, Wojciech Kryszewski, Jacobo Pejsachowicz and Robert Skiba for the generous support with information.

$$ \begin{array}{lcccr}Department \, of \, Mechanics \, and \,  Mathematics & & & & Department \, of \, Mathematics \,and \,  Informatics \\
Moscow \, State \,  University & & & &  University \, of \, Bremen \\
Leninskie \, Gory   & & & &         Universitaetsallee \, 29 \\
119992 \, Moscow  & & & &       28359 \, Bremen \\
Russia      & & & &         Germany \\
egskl@higeom.math.msu.su & & & & skordev@cevis.uni-bremen.de
\end{array}$$
\end{document}